\documentclass[11pt]{amsart}
\usepackage[T1]{fontenc}
\usepackage{lmodern,microtype,amsmath,amssymb,mathtools,enumitem}
\usepackage{tikz-cd}
\usepackage[a4paper,margin=26mm]{geometry}
\usepackage[colorlinks=true,linkcolor=blue,citecolor=blue,urlcolor=blue]{hyperref}
\newtheorem{theorem}{Theorem}[section]
\newtheorem{proposition}[theorem]{Proposition}
\newtheorem{lemma}[theorem]{Lemma}
\newtheorem{corollary}[theorem]{Corollary}
\theoremstyle{definition}
\newtheorem{definition}[theorem]{Definition}
\theoremstyle{remark}

\newcommand{\PP}{\mathbb P}
\newcommand{\Aff}{\mathbb A}
\newcommand{\CC}{\mathbb C}
\newcommand{\QQ}{\mathbb Q}

\newcommand{\NN}{\mathbb N}
\newcommand{\Gm}{\mathbb G_m}
\newcommand{\OO}{\mathcal O}
\newcommand{\cS}{\mathcal S}
\newcommand{\cE}{\mathcal E}

\newcommand{\cU}{\mathcal U}
\newcommand{\cD}{\mathcal D}
\newcommand{\cC}{\mathcal C}
\newcommand{\reg}{\mathrm{reg}}

\newcommand{\Supp}{\operatorname{Supp}}
\newcommand{\Spec}{\operatorname{Spec}}

\newcommand{\Hom}{\operatorname{Hom}}
\newcommand{\Diff}{\operatorname{Diff}}

\newcommand{\rk}{\operatorname{rank}}
\newcommand{\lcm}{\operatorname{lcm}}

\numberwithin{equation}{section}
\title{Free rational curves on singular Fano varieties}
\hypersetup{pdftitle={Free rational curves on singular Fano varieties},pdfauthor={Zheng Xu}}
\author{Zheng Xu}
\address{Department of Mathematics, Southern University of Science and Technology (SUSTech), 1088 Xueyuan Rd, Nanshan District, Shenzhen, Guangdong, China, 518055}
\email{xuz6@sustech.edu.cn}
\date{}
\begin{document}
\begin{abstract}
We prove that every klt Fano variety admits a free rational curve in its smooth locus.
\end{abstract}
\maketitle

\section{Introduction}\label{sec:introduction}

For a smooth variety $V$,
a morphism $f:\PP^1\to V$ is
\emph{free} if $f$ is nonconstant and $f^*T_V$ is globally generated.
In this paper, we prove the following folklore conjecture.

\begin{theorem}\label{mainthm}
Let $X$ be a klt Fano variety. There exists a free morphism
\[
 f:\PP^1\longrightarrow X_{\reg}.
\]
\end{theorem}

It is known to experts that Theorem \ref{mainthm} implies the existence of very free rational curves. Then the following result follows.

\begin{corollary}
    Let $X$ be a klt Fano variety. Then $X_{\reg}$ is rationally connected.
\end{corollary}

Keel and McKernan proved that the smooth locus of a log Fano surface
is rationally connected \cite{KM97}. In dimension three,
Jovinelly, Lehmann and Riedl constructed a very free rational curve in
$X_{\reg}$ for every terminal Fano threefold $X$ over an algebraically
closed field of characteristic zero \cite{JLR25,JLR26}.

Our proof uses logarithmic bend-and-break,
conjectured by Keel and McKernan \cite{KM97}.
For a smooth complex projective variety $V$ with reduced SNC divisor
$D$, Chen and Tseng proved that through a general point of an integral
curve $C\not\subset D$ with $(K_V+D)\cdot C<0$ there is a rational
curve whose normalization meets $D$ at most once
\cite{CT26a}.
Their proof is based on the weighted jet construction of Guo and Xie
\cite{GX26}. Our main new idea is to use the minimal model program (MMP) to descend free rational curves from a log resolution. To do this, one has to establish an inductive statement about dlt weak log Fano pairs (see Theorem \ref{thm:induction}).

\subsection*{Outline of the proof}

Let $(X,\Gamma)$ be a dlt weak log Fano pair as in
Theorem~\ref{thm:induction}, and fix a very ample divisor $H$ on $X$.
We descend a rational curve from a log resolution of $X$ to a small $\QQ$-factorization of $Y_N\to X$ along the relative minimal model program of Lemma~\ref{lem:mmp}.
A general free curve avoids a flipping locus because it has codimension
at least two. At a divisorial contraction $\phi:(Y,\Gamma_Y)\to Y'$
with exceptional prime divisor $E$,
a curve of positive $E$-degree meets $E$ in every deformation.

Let $p:\cS\to Y$ be the rooted model and $\rho:\cS\to X$ the
composite with $Y\to X$. Set $\cE=(p^{-1}E)_{\rm red}$ and
$D'=(p^{-1}\lfloor\Gamma_Y\rfloor)_{\rm red}-\cE$. Deform the
logarithmically free curve $f:\PP^1\to(\cS,\cE+D')$ so that
$f^*\cE=\sum_{a=1}^r[p_a]$ with distinct $p_a$ and $f(p_a)\notin D'$.
At each $p_a$ we attach a logarithmically free curve
$h_a:\PP^1\to(\cE,D'|_{\cE})$ through $f(p_a)$ with $\rho h_a$
constant and $h_a^*N_{\cE/\cS}\simeq\OO(-1)$.
Proposition~\ref{prop:cancellation} smooths the nodal map to a curve
$g$ logarithmically free along $D'$. Then we have
\[
 \deg g^*\OO_{\cS}(\cE)=r-r=0,\qquad
 \deg(\rho g)^*\OO_X(H)=\deg(\rho f)^*\OO_X(H)>0.
\]
A marked point of $g$ lies outside $\cE$, so $g^*\cE$ is effective
of degree zero. Thus $g$ doesn't meet $\cE$ and passes to the next model
with positive degree over $X$.

The curves $h_a$ come from induction on dimension. Write
$E\to W\to\phi(E)$ for the Stein factorization of $\phi|_E$.
For a general fibre $F$ of $E\to W$, adjunction gives a dlt log
Fano pair $(F,\Gamma_F)$. The inductive statement allows arbitrary root
orders along $\lfloor\Gamma_F\rfloor$. Let $p^0:\cS^0\to Y$ be
the model with all prescribed roots except the one along $E$, and set
$\cE^0=((p^0)^{-1}E)_{\rm red}$. Lemma~\ref{lem:fibre} gives
vertical curves through a dense open subset of $\cE^0$, all of the same
normal degree $-d<0$. On the $d$-th root along $\cE^0$, each curve
lifts with normal line $\OO(-1)$ because
$\OO(-1)^{\otimes d}\simeq\OO(-d)$
(Proposition~\ref{prop:root-orders}). We fix the root orders of the other
exceptional prime divisors on $Y$ first, since the MMP contracts them later.

With all orders fixed, take general complete intersection curves of
class $H^{n-1}$, where $n=\dim X$. Their strict transforms lift to the
rooted resolution after finite covers. Since $(K_X+\Gamma)\cdot H^{n-1}<0$,
logarithmic bend-and-break in Lemma~\ref{lem:initial} gives the initial
logarithmically free curve of positive degree over $X$. We then follow
the flips and contractions in forward order. For $\Gamma=0$ the small
morphism $Y_N\to X$ is an isomorphism over $X_{\reg}$, and its inverse
image of $X\setminus X_{\reg}$ has codimension at least two. A general
deformation therefore gives a free map $\PP^1\to X_{\reg}$.

\subsection*{Acknowledgements}
I would like to thank Zhiyu Tian for helpful discussions, and Wenhao Ou and Chenyang Xu for encouragement.
The author is supported by start-up and matching funds (Y01286160) from the Southern University of Science and Technology (SUSTech).

\subsection*{AI disclosure}

In April 2026, I began using ideas of the minimal model program, adjunction and induction on dimension, the comb argument, stacks, and boundedness to study Theorem \ref{mainthm}. I obtained a much weaker version of \cite[Theorem 1.1]{CT26a} in July 2026. After the appearance of \cite{CT26a}, I realized that \cite[Theorem 1.1]{CT26a}, combined with my previous ideas, could imply Theorem \ref{mainthm}. I then asked GPT 6 Astra to generate a proof with these ideas, which turned out to be more complicated than I had imagined. I checked, simplified, and rewrote the proof. During this process, I also used GPT 6.2 Astra to review the results and polish the English.

\section{Preliminaries}\label{sec:preliminaries}

We work over $\CC$. We refer to \cite{KM98}, \cite{Kol13} and \cite{BCHM08} for basic knowledge of birational geometry. In this paper, we use $\eta$ to denote the generic point.

\subsection{Log discrepancies and adjunction}

For a normal pair $(Y,\Gamma)$, $K_Y+\Gamma$ is $\QQ$-Cartier, and $\Gamma$ is an effective $\QQ$-divisor with coefficients in $[0,1]$. On a
birational model $\mu:Z\to Y$ we write
\[
 A_{Y,\Gamma}(F)=1+\operatorname{coeff}_F(K_Z-\mu^*(K_Y+\Gamma))
\]
for the log discrepancy of a prime divisor $F$.

\begin{proposition}\label{prop:log discrepancies under finite covers} (\cite[Proposition~5.20]{KM98})
Let $\nu:Z\to Y$ be a finite surjective morphism of normal varieties.
Assume that $K_Z+\Gamma_Z=\nu^*(K_Y+\Gamma_Y)$. For each prime
divisor $F$ over $Z$, write $\operatorname{ord}_F|_{K(Y)}=e\operatorname{ord}_P$,
where $P$ is a prime divisor over $Y$ and $e$ is the ramification index.
Let $c_Z(F)$ and $c_Y(P)$ denote the respective centres. Then
\[
 A_{Z,\Gamma_Z}(F)=e\,A_{Y,\Gamma_Y}(P),\qquad \nu(c_Z(F))=c_Y(P).
\]
In particular, $(Z,\Gamma_Z)$ is klt (respectively plt, lc) if and
only if $(Y,\Gamma_Y)$ is. If $(Z,\Gamma_Z)$ is lc, then for each
lc centre $W$ of $(Z,\Gamma_Z)$ the image $\nu(W)$ is an lc centre
of $(Y,\Gamma_Y)$.
\end{proposition}

\begin{proposition}\label{prop:adjunction} (\cite[Theorems~4.9 and~4.16, (4.2.7)]{Kol13})
Let $(Y,S+\Lambda)$ be a dlt pair with $S$ a prime divisor of
coefficient one. Then $S$ is normal. For sufficiently divisible $m$,
the $m$-th power of the Poincar\'e residue along $S$ induces the
adjunction formula
\[
 (K_Y+S+\Lambda)|_S=K_S+\Diff_S(\Lambda).
\]
The pair $(S,\Diff_S(\Lambda))$ is dlt. If $(Y,S+\Lambda)$ is plt
near $S$, then $(S,\Diff_S(\Lambda))$ is klt. Let $W$ be an
irreducible component of an intersection of distinct coefficient-one
divisors of $S+\Lambda$. Then $W$ is an lc centre of $(Y,S+\Lambda)$,
and $(Y,S+\Lambda)$ is SNC at $\eta_W$.
\end{proposition}

\subsection{Root stacks and logarithmic tangent bundles}

\begin{definition}\label{def:root-stack} (\cite[Example~2.7]{Cad07})
Let $V$ be a Deligne--Mumford stack, let $L$ be a line
bundle on $V$, let $s\in H^0(V,L)$, and let $r\geq 1$ be an integer.
The \emph{$r$-th root stack} of $(L,s)$ is the stack
\[
 p:\mathcal R=\sqrt[r]{(L,s)/V}\longrightarrow V
\]
defined as follows. For a scheme $T$ and a morphism $u:T\to V$, an
object over $u$ is a triple $(M,t,\varphi)$ consisting of a line
bundle $M$ on $T$, a section $t\in H^0(T,M)$, and an isomorphism of
line bundles on $T$
\[
 \varphi:M^{\otimes r}\xrightarrow{\sim}u^*L
 \qquad\text{such that}\qquad
 \varphi(t^{\otimes r})=u^*s.
\]
An isomorphism from $(M,t,\varphi)$ to $(M',t',\varphi')$ over the
same $u$ is an isomorphism $\alpha:M\to M'$ on $T$ satisfying
\[
 \alpha(t)=t',
 \qquad
 \varphi'\circ\alpha^{\otimes r}=\varphi.
\]
Pullback along $a:T'\to T$ sends this triple to
$(a^*M,a^*t,a^*\varphi)$ over $u\circ a$. The universal root line
bundle $\mathcal L$ and its section $\tau$ are defined on $\mathcal R$.
They satisfy
\[
 (\mathcal L^{\otimes r},\tau^{\otimes r})\simeq(p^*L,p^*s).
\]
For an effective Cartier divisor $A$ on $V$, let $s_A$ be the image
of $1$ under $\OO_V\hookrightarrow\OO_V(A)$. We write
$\sqrt[r]{(V,A)}$ for the root stack of $(\OO_V(A),s_A)$.
For labelled effective Cartier divisors $A_1,\ldots,A_k$ on $V$
with prescribed orders $r_1,\ldots,r_k$, their \emph{separate root
stack} is the fibre product
\[
 \sqrt[r_1]{(V,A_1)}\times_V\cdots\times_V\sqrt[r_k]{(V,A_k)}.
\]
A root of order one is canonically isomorphic to $V$ over $V$.
\end{definition}

A rational curve on a stack $V$ is a morphism $f:\PP^1\to V$ whose
composite with the coarse moduli morphism is nonconstant. If $V$ is
smooth, such a curve is \emph{free} if $f^*T_V$ is globally generated
on $\PP^1$.

\begin{definition}\label{def:logfree}
For a smooth stack $V$ with reduced SNC divisor
$D$, a morphism $f:\PP^1\to V$ is \emph{logarithmically free} if
its coarse map is nonconstant, its image is not contained in $D$, and
$f^*T_V(-\log D)$ is globally generated. 
\end{definition}

\subsection{Deformations of maps from curves}

\begin{lemma}\label{lem:p1}
Let $E$ be a vector bundle on $\PP^1$. Write
$\operatorname{ev}_p:H^0(\PP^1,E)\to E|_p$ for the evaluation map at $p$. The following are equivalent.

\noindent (1) $\operatorname{ev}_p$ is surjective for some $p$.

\noindent (2) $E$ is globally generated.

\noindent (3) $H^1(\PP^1,E(-1))=0$.

\noindent If \textup{(1)}--\textup{(3)} hold, then $H^1(\PP^1,E)=0$.
\end{lemma}
\begin{proof}
For $q\in\PP^1$ we have $\OO_{\PP^1}(-q)\simeq\OO_{\PP^1}(-1)$. The exact sequence
$0\to E(-q)\to E\to E|_q\to0$ shows that \textup{(3)} implies
surjectivity of $\operatorname{ev}_q$ and $H^1(\PP^1,E)=0$.
Thus \textup{(3)} implies \textup{(2)}, and \textup{(2)} implies \textup{(1)}.
Assume that $\operatorname{ev}_p$ is surjective for some $p\in\PP^1$.
Take $r=\rk E$ global sections whose
values form a basis of $E|_p$. They give an exact sequence
\[
 0\longrightarrow\OO_{\PP^1}^{r}\longrightarrow E\longrightarrow Q
 \longrightarrow0,
\]
where $Q$ is supported on finitely many points. Tensoring by $\OO_{\PP^1}(-1)$
gives $H^1(\PP^1,E(-1))=0$, since both $H^1(\PP^1,\OO_{\PP^1}(-1))$
and $H^1(\PP^1,Q(-1))$ vanish. This proves \textup{(3)}.
\end{proof}

\begin{proposition}\label{prop:sections} (\cite{Ols06})
Let $B$ be a scheme of finite type over $\CC$. Let $\pi:C\to B$ be a
proper flat finitely presented family of nodal curves, with $C$ a
scheme. Let $V$ be a smooth separated Deligne--Mumford stack of
finite type. For
$b\in B(\CC)$ and $h:C_b\to V$, we write $E=h^*T_V$.
Write $M=\underline{\Hom}_B(C,V\times B)$. 
Then $M$ is a Deligne--Mumford
stack locally of finite presentation over $B$. Moreover, the following hold.

\noindent (1) If $H^1(C_b,E)=0$, then $M\to B$ is smooth at $h$. After an
\'etale base change near $b$, the map $h$ extends to the family.

\noindent (2) If a section $z:B\to C$ of $\pi$ meets the relative smooth locus and
$H^1(C_b,E(-z_b))=0$, then the relative evaluation
$(\operatorname{ev}_z,\pi_M):M\to V\times B$ is smooth at $h$, where
$\pi_M:M\to B$ is the structure map.
\end{proposition}

We next apply ordinary evaluation to move a logarithmically free
curve away from a closed subset and to arrange transverse contacts.

\begin{lemma}\label{lem:free}
Let $(V,D)$ be a smooth Deligne--Mumford stack with reduced simple
normal crossing boundary. Assume that $f:\PP^1\to(V,D)$ is a logarithmically
free morphism. Then the following hold.

\noindent (1) $f$ is free as an ordinary map. 

\noindent (2) A general ordinary deformation
of $f$ remains logarithmically free, avoids any fixed closed subset of
codimension at least two, and meets each component of $D$ transversely
at distinct points away from the other components.
\end{lemma}
\begin{proof}

\noindent Take $p\in\PP^1\setminus f^{-1}D$. The inclusion
$f^*T_V(-\log D)\to f^*T_V$ is an isomorphism at $p$. Since
$f^*T_V(-\log D)$ is globally generated, the global sections of
$f^*T_V$ span its fibre at $p$.
Lemma~\ref{lem:p1} gives $H^1(f^*T_V(-q))=0$ for every $q\in\PP^1$. This proves (1).

By Proposition~\ref{prop:sections}\textup{(2)} with $B=\Spec\CC$,
evaluation at each $q\in\PP^1$ is smooth at $f$. Take a smooth
finite-type parameter atlas $T$ with a point $t_0$ representing $f$.
The universal evaluation $e:\PP^1\times T\to V$ is smooth along
$\PP^1\times\{t_0\}$. Since $\PP^1\times T\to T$ is proper, we may
shrink $T$ about $t_0$ so that $e$ is smooth.
Write $f_t=e|_{\PP^1\times\{t\}}$ and let $\bar f_t$ denote its
coarse map. We may assume that $T$ is irreducible. Take $q\in\PP^1$ such that
$\bar f(p)\ne\bar f(q)$. The conditions
\[
 f_t(p)\notin D,\qquad
 \bar f_t(p)\ne\bar f_t(q),\qquad
 H^1\bigl(\PP^1,f_t^*T_V(-\log D)\otimes\OO_{\PP^1}(-1)\bigr)=0
\]
hold at $t_0$ by the choice of $p,q$ and Lemma~\ref{lem:p1}.
They are open by closedness of $D$, separatedness of the coarse
space, and semicontinuity, respectively. Shrink $T$ to this open
neighbourhood. By Lemma~\ref{lem:p1}, every $f_t$ is logarithmically
free.

Fix a closed subset $Z\subset V$ of codimension at least two.
Write $D=\sum_i D_i$ and set
\[
 W=Z\cup\bigcup_{i<j}(D_i\cap D_j).
\]
Since $\operatorname{codim}_V W\geq2$ and $e$ is smooth,
$\dim e^{-1}(W)\leq\dim T-1$. The projection
$\pi:\PP^1\times T\to T$ is proper, so
$U=T\setminus\pi(e^{-1}(W))$ is dense and open. For every $t\in U$,
the curve $f_t$ avoids $W$.

Over $U$, the scheme-theoretic inverse image
$Y=(e|_{\PP^1\times U})^{-1}(D)$ is smooth because $e$ is smooth
and $D\setminus W$ is smooth. Since $f_t(p)\notin D$, the projection $Y\to U$ is proper
with finite fibres, and hence is finite.
By generic smoothness, it is finite \'etale
over a dense open subset of $U$. For every geometric point $t$ of
this open, the divisor $f_t^*D$ is therefore reduced. Each point
$x$ of its support maps to a unique $D_i$ and satisfies
$\operatorname{ord}_x(f_t^*D_i)=1$. Thus all contacts are transverse
and occur at distinct points of $\PP^1$. This concludes (2).
\end{proof}

\subsection{Results of Chen and Tseng on rational curves}\label{subsec:chen-tseng}

\begin{theorem}\label{thm:ct-contact} (\cite[Theorem~1.1]{CT26a})
Let $Y$ be a smooth and projective variety, $B$ a reduced SNC divisor on
$Y$, and $C\subset Y$ an integral curve meeting
$Y\setminus B$. If $(K_Y+B)\cdot C<0$, there is a dense open set
$C^\circ\subset C\setminus B$ such that for every $y\in C^\circ(\CC)$
there is a morphism $\nu:\PP^1\to Y$, birational onto its image, such that
\[
 \nu(0)=y,\qquad \Supp(\nu^*B)\subset\{\infty\}.
\]
\end{theorem}

For the root stack constructions, let $A_1+\cdots+A_l$ be a reduced
SNC divisor on a smooth variety $V$. Write $p:\cS\to V$ for the
separate roots of positive orders $n_i$, $R_i=(p^{-1}A_i)_{\rm red}$, and
$V^*=V\setminus\bigcup_i A_i$. Then for $J\subset\{1,\ldots,l\}$ we have
\begin{equation}\label{eq:root-canonical}
 K_{\cS}+\sum_{i\in J}R_i
 =p^*\left(K_V+\sum_{i\notin J}(1-1/n_i)A_i+\sum_{i\in J}A_i\right).
\end{equation}

\begin{lemma}\label{lem:ct-lift} (\cite[\S3.1]{CT26b})
Let $F:C\to V$ be a map from a smooth connected curve whose image meets
$V^*$. Then $F$ lifts to $\cS$ if and only if
 $n_i\mid \operatorname{ord}_c(F^*A_i)$ for all $i$ and $c\in C(\CC)$.
 
 If $C$ is projective, some finite surjective map $\tau:C'\to C$
from a smooth connected projective curve makes $F\tau$ lift to $\cS$.
\end{lemma}

To apply Theorem~\ref{thm:ct-contact} to a map $b:C\to\cS$, we lift
$b$ after a finite cover of $C$ to a smooth quasi-projective scheme.

\begin{proposition}\label{prop:ct-torsor} (\cite[Lemma~3.1 and proof of Proposition~3.2]{CT26b})
Let $M_1,\ldots,M_l$ be line bundles on a smooth connected projective
curve $C$. There are a finite surjection $\lambda:C'\to C$ from a
smooth connected projective curve, a map $a:C'\to A$ to an abelian
variety, and line bundles $L_i$ on $A$ with $\lambda^*M_i\simeq a^*L_i$.

Assume that $V$ is projective and $M_i=b^*\OO_{\cS}(R_i)$ for a map
$b:C\to\cS$. On $P:=A\times\cS$, with projections $\pi_A,\pi_{\cS}$, we set
\[
 W_i=\operatorname{Isom}_{P}(\pi_A^*L_i,\pi_{\cS}^*\OO_{\cS}(R_i)),
 \qquad W=W_1\times_P\cdots\times_P W_l .
\]
Then $q:W\to P$ is a $\Gm^l$-torsor and $W$ is a smooth
quasi-projective scheme. The isomorphisms $\lambda^*M_i\simeq a^*L_i$ give
$\widetilde b:C'\to W$ with $q\widetilde b=(a,b\lambda)$.
For $l=0$ take $C'=C$, $A$ a point, and $W=V$.
\end{proposition}

We will use the following families of maps.

\begin{lemma}\label{lem:ct-families}
Assume that $V$ is projective. There are countably many finite-type
$\CC$-schemes $T_e$ and morphisms $\Phi_e:\PP^1\times T_e\to\cS$
with the following property. For every ordinary map $f:\PP^1\to\cS$
with nonconstant coarse map $pf$, there are $e\geq1$ and $t\in T_e(\CC)$
such that $\Phi_e|_{\PP^1\times\{t\}}\simeq f$.
\end{lemma}
\begin{proof}
Fix a very ample line bundle $H$ on $V$. We record the parameter of
a map $f$ of coarse $H$-degree $e>0$ by taking its graph
\[
 (\mathrm{id},f):\PP^1\longrightarrow\PP^1\times\cS
\]
of degree $1+e$ for $\OO_{\PP^1}(1)\boxtimes H$ on $\PP^1\times V$.
Any automorphism of the coarse graph fixes its first projection,
which is the identity. Thus the graph is an unmarked stable map.

Since $\PP^1\times\cS$ is proper and tame with projective coarse
space, the stack of genus-zero twisted stable maps of degree $1+e$
is of finite type. An unmarked twisted curve has stabilizers only at
nodes \cite[Theorem~1.4.1 and Definition~4.1.2]{AV02}. Hence over a
finite-type scheme atlas $U_e$ of the smooth-domain locus the universal
curve $\pi:C_e\to U_e$ is an ordinary smooth proper family.
Write $(a,b):C_e\to\PP^1\times\cS$ for the universal map.

Let $T_e\subset U_e$ be the locus of points $t$ for which
$a_t:C_{e,t}\to\PP^1_{\kappa(t)}$ is an isomorphism. Since $C_e$ and
$\PP^1\times U_e$ are smooth and proper over $U_e$, $T_e$ is open
and of finite type, and
$(a,\pi):C_e|_{T_e}\xrightarrow{\sim}\PP^1\times T_e$
by \cite[Tag~05XD]{Sta26}. Set $\Phi_e=b\circ(a,\pi)^{-1}$.
For $t\in T_e(\CC)$ the map $p\Phi_{e,t}$ has $H$-degree $e$,
since its graph has degree $1+e$.
The graph of any $f$ of coarse degree $e$ is represented
by some $t\in U_e(\CC)$ and an isomorphism $i:\PP^1\to C_{e,t}$
with $a_t i=\mathrm{id}$ and $b_t i\simeq f$. Hence $t\in T_e$ and
$\Phi_{e,t}=b_t a_t^{-1}\simeq f$. There is one such family for each
positive integer $e$.
\end{proof}

\section{The inductive statement and the initial curve}\label{sec:inductive-pair}

In this section, we give the inductive statement. We construct a logarithmically free curve, with prescribed root orders on the boundary, on a root stack over a smooth projective variety. We then apply the construction to a curve to solve the case $\dim X=1$.

\subsection{The inductive statement}

Let $I_{\rm st}=\{1-1/r:r\in\NN_{>0}\}\cup\{1\}$ be the standard coefficient set. For a dlt pair
$(X,\Gamma)$ with coefficients in $I_{\rm st}$, we write
\begin{equation}\label{eq:split}
 \Gamma=\Delta+B,
 \qquad B=\lfloor\Gamma\rfloor=\sum_{j=1}^{b}B_j,
 \qquad
 \Delta=\sum_{k=1}^{a}(1-1/r_k)\Delta_k,
 \quad r_k\geq2.
\end{equation}
All $B_j,\Delta_k$ are distinct prime divisors. Let $U$ be the open
locus where $X$ is smooth and $\Supp\Gamma$ is SNC. By the definition of a dlt pair,
$\operatorname{codim}_X(X\setminus U)\geq2$. For positive
integers $m_1,\dots,m_b$ let $p:\cU\to U$ be the projection from
\begin{equation}\label{eq:root-open}
 \cU=
 \sqrt[r_1]{(U,\Delta_1|_U)}\mathbin{\times_U}\cdots
 \mathbin{\times_U}\sqrt[r_a]{(U,\Delta_a|_U)}
 \mathbin{\times_U}\prod_{j=1}^{b}\sqrt[m_j]{(U,B_j|_U)}.
\end{equation}
Set $B_{\cU}=(p^{-1}(B|_U))_{\rm red}$. For $n\geq1$ denote the
following statement by $Q_n$.

\begin{theorem}\label{thm:induction}
Let $(X,\Delta+B)$ be an $n$-dimensional projective dlt pair as in
\eqref{eq:split}. Assume that
\[
 -(K_X+\Delta+B)\text{ is nef and big},
 \qquad B\text{ is }\QQ\text{-Cartier}.
\]
For every choice of the integers $m_j$ in \eqref{eq:root-open},
there is a logarithmically free ordinary morphism
\[
 f:\PP^1\longrightarrow (\cU,B_{\cU}).
\]
Its ordinary deformations cover a dense open subset and can avoid any
specified additional closed subset of codimension at least two in
$\cU$.
\end{theorem}

 For \eqref{eq:root-open} the canonical divisor formula
\eqref{eq:root-canonical} reads
$K_{\cU}+B_{\cU}=p^*((K_X+\Delta+B)|_U)$.

\subsection{A curve on a smooth root stack}\label{sec:initial}

Let $p:\cS\to V$ be the separate root stack of a smooth projective
variety along an SNC divisor $A_1+\cdots+A_l$, with orders $n_i$.
Write $R_i=(p^{-1}A_i)_{\rm red}$, fix $D=\sum_{i\in J}R_i$, and
set $V^*=V\setminus\bigcup_i A_i$. Chen and Tseng prove the case
$D=0$ in \cite[Proposition~3.2]{CT26b}. To control contact with $D$ we use
logarithmic bend-and-break on the torus torsor of
Proposition~\ref{prop:ct-torsor}.

\begin{lemma}\label{lem:one-contact}
Let $b:C\to\cS$ be a map from a smooth connected projective curve
with $pb(C)\cap V^*\ne\varnothing$ and $\deg b^*(K_{\cS}+D)<0$.
For a general $c\in C$ there is an ordinary map
$f:\PP^1\to\cS$ with $f(0)\simeq b(c)$, $pf$ nonconstant and
$f^{-1}D\subset\{\infty\}$.
\end{lemma}
\begin{proof}
By Proposition~\ref{prop:ct-torsor}, take a finite cover $\lambda:C'\to C$
and a lift $\widetilde b:C'\to W$ with $h\widetilde b\simeq b\lambda$.
Here $h=\pi_{\cS}\circ q$ for the torus torsor $q:W\to A\times\cS$.
We construct $g$ and $f$ in the following diagram, with
$F|_{\Aff^1}=phg$ and $pf=F\rho$.
\[
\begin{tikzcd}[row sep=2em,column sep=3.3em]
 C' \arrow[r,"\widetilde b"] \arrow[d,"\lambda"'] &
 W \arrow[d,"h"] &
 \Aff^1 \arrow[l,"g"'] \arrow[r,hook] &
 \PP^1 \arrow[d,"F"] &
 \PP^1 \arrow[l,"\rho"'] \arrow[d,dashed,"f"] \\
 C \arrow[r,"b"'] & \cS \arrow[rr,"p"'] & & V & \cS \arrow[l,"p"]
\end{tikzcd}
\]
By embedded resolution, take a smooth projective compactification $\overline W$ such that
$\overline D=(\overline W\setminus W)\cup\overline{h^{-1}D}$ is a
reduced SNC divisor. Since $\omega_{W/(A\times\cS)}$ and $\omega_A$
are trivial, $\omega_{\overline W}(\overline D)|_W\simeq h^*\omega_{\cS}(D)$.
The projection formula gives
\[
 (K_{\overline W}+\overline D)\cdot\widetilde b_*[C']
    =\deg(\lambda)\deg b^*(K_{\cS}+D)<0 .
\]
Thus the reduced image of $\widetilde b$ is a curve of negative
logarithmic degree. Since $pb(C)$ meets $V^*$, the point
$pb(\lambda(c'))$ lies in $V^*$ for general $c'\in C'$. In particular
the image curve meets $\overline W\setminus\overline D$. Apply
Theorem~\ref{thm:ct-contact} to this curve and restrict the resulting
map to $\Aff^1$. For general $c'$ we obtain a nonconstant map
$g:\Aff^1\to W\setminus h^{-1}D$ with $g(0)=\widetilde b(c')$.

Since $V$ is proper, $phg$ extends to $F:\PP^1\to V$.
If $F$ were constant, its value would lie in $V^*$, where $p$ is an
isomorphism. The map $\pi_A qg$ is also constant, since it extends to
$\PP^1$ and $A$ is an abelian variety. Then $g$ would map into a
single fibre of $q$, isomorphic to $\Gm^l$. Every map
$\Aff^1\to\Gm^l$ is constant because $\CC[t]^\times=\CC^\times$.
Hence $F$ is nonconstant.

On $\Aff^1$ the lift $hg$ gives $n_i\mid\operatorname{ord}_z F^*A_i$
for every $i$ and every $z\in\Aff^1$. Take $N=\lcm_i n_i$, with $N=1$
if $l=0$, and set $\rho(t)=t^N$. Pullback by $\rho$ preserves these
divisibilities on $\Aff^1$. At $\infty$ we have
\[
 \operatorname{ord}_\infty(F\rho)^*A_i=N\operatorname{ord}_\infty F^*A_i .
\]
Lemma~\ref{lem:ct-lift} gives $f:\PP^1\to\cS$ with $pf=F\rho$.
This coarse map is nonconstant. Since $\rho(0)=0$ and $p$ is an
isomorphism over $V^*$, $f(0)\simeq b(\lambda(c'))$. Moreover,
$F(\Aff^1)$ misses $p(D)$, so $f^{-1}D\subset\{\infty\}$.
The finite set excluded from $C'$ has finite image in $C$, and
$\lambda$ is surjective. This proves the assertion for general $c\in C$.
\end{proof}

\begin{lemma}\label{lem:initial}
Assume that through a general point of $V$
there is a projective curve $b:C\to\cS$ with
$\deg b^*(K_{\cS}+D)<0$. Then there exists a logarithmically free
ordinary rational curve $f:\PP^1\to(\cS,D)$ with
$f^{-1}D\subset\{\infty\}$.
\end{lemma}
\begin{proof}
For each family $\Phi:\PP^1\times T\to\cS$ in Lemma~\ref{lem:ct-families},
restrict $T$ to the open set where $p\Phi(0,t)\in V^*$.
The condition $\Phi_t^{-1}D\subset\{\infty\}$ defines the constructible set
\[
 T\setminus\operatorname{pr}_T\bigl((\Aff^1\times T)\cap\Phi^{-1}D\bigr)
\]
by Chevalley's theorem \cite[Tag~054K]{Sta26}. Write each such set as a
finite union of smooth integral locally closed schemes $T_\nu$. Denote
the restricted families by $\Phi_\nu$ and set $e_\nu(t)=p\Phi_\nu(0,t)$.

Let $U_0\subset V^*$ be a dense open subset such that the hypothesis
holds at every closed point of $U_0$. We claim that some
$e_\nu:T_\nu\to V$ is dominant. Otherwise the
closures $Z_\nu=\overline{e_\nu(T_\nu)}$ are countably many proper
closed subsets of $V$. Since $\CC$ is uncountable, we can take a
point $x\in U_0(\CC)\setminus\bigcup_\nu Z_\nu$.
Take $b:C\to\cS$ as in the assumption with $x\in pb(C)$, replacing
$C$ by its normalization. Each $(pb)^{-1}Z_\nu$ is finite
since $x\in pb(C)\setminus Z_\nu$. By Lemma~\ref{lem:one-contact},
for general $c\in C$ there is a map in one of these families with
$pb(c)=e_\nu(t)$. A dense open subset of $C(\CC)$ would therefore lie
in a countable union of finite sets. This proves the claim.

Take $T_\nu$ with $e_\nu$ dominant. Fix $t\in T_\nu(\CC)$ and set
$f=\Phi_{\nu,t}$. For $v\in T_tT_\nu$, the differential of $\Phi_\nu$
gives the section $s_v(q)=d\Phi_{\nu,(q,t)}(0,v)$ of $f^*T_{\cS}$.
For $i\in J$ the effective Cartier divisor $\Phi_\nu^*R_i$ is supported
on the prime divisor $\{\infty\}\times T_\nu$. Take a parameter $z$ at
$\infty$ and an \'etale local equation $x_i=0$ for each $R_i$ through
$f(\infty)$. Near $(\infty,t)$ we have $\Phi_\nu^*x_i=u_i z^{k_i}$,
where $u_i$ is a unit and $k_i\geq0$. Since $\partial_v z=0$, the
$x_i$-component of $s_v$ is
\[
 \left.\partial_v(\Phi_\nu^*x_i)\right|_t
     =\left.\frac{\partial_v u_i}{u_i}\right|_t f^*x_i .
\]
Hence $s_v$ is logarithmically tangent to each $R_i$ at $\infty$.
Since $f(\Aff^1)\cap D=\varnothing$, we obtain the factorization
\[
 T_tT_\nu \xrightarrow{\,v\mapsto s_v\,} H^0(\PP^1,f^*T_{\cS}(-\log D))
 \xrightarrow{\,\operatorname{ev}_0\,} T_{f(0)}\cS
 \xrightarrow{\,dp_{f(0)}\,} T_{e_\nu(t)}V
\]
of $de_{\nu,t}$, since $e_\nu(t)=p(\Phi_\nu(0,t))$. The map
$dp_{f(0)}$ is an isomorphism because $p(f(0))\in V^*$.
For general $t$, generic smoothness gives surjectivity of $de_{\nu,t}$,
so $\operatorname{ev}_0$ is surjective. By Lemma~\ref{lem:p1},
$f^*T_{\cS}(-\log D)$ is globally generated. The coarse map $pf$ has
positive degree by Lemma~\ref{lem:ct-families}, and $f(0)\notin D$.
\end{proof}

To solve the one-dimensional case, we apply Lemma~\ref{lem:initial}
to a root stack over a curve.

\begin{proposition}\label{prop:curve}
Let $X=C$ be a normal projective curve with $\deg(K_C+\Delta+B)<0$.
For every choice of the $m_j$ there exists $f:\PP^1\to\cU$ of
positive coarse degree $e$ with
\begin{equation}\label{eq:tangent-degree}
 f^*T_{\cU}(-\log B_{\cU})\simeq\OO_{\PP^1}(-e\deg(K_C+\Delta+B)).
\end{equation}
\end{proposition}
\begin{proof}
Since $C$ is a normal curve, it is smooth and $U=C$. Apply
Lemma~\ref{lem:ct-lift} to $\mathrm{id}_C$. We obtain a finite cover
$\tau:C'\to C$ and a map $b:C'\to\cU$ with $pb=\tau$.
By \eqref{eq:root-canonical},
\[
 \deg b^*(K_{\cU}+B_{\cU})=(\deg\tau)\deg(K_C+\Delta+B)<0 .
\]
Since $\tau$ is surjective, Lemma~\ref{lem:initial} gives a
logarithmically free map $f:\PP^1\to(\cU,B_{\cU})$. Set $e=\deg(pf)>0$.
The logarithmic tangent line is $\omega_{\cU}(B_{\cU})^\vee$. Its
pullback has degree $-e\deg(K_C+\Delta+B)$ by \eqref{eq:root-canonical}.
A line bundle on $\PP^1$ is determined by its degree, which proves
\eqref{eq:tangent-degree}.
\end{proof}

Since $pf$ has positive degree, it is surjective. A closed subset of
codimension at least two in the curve $\cU$ is empty. This proves $Q_1$.

\section{Birational models and curves on exceptional fibres}\label{sec:models}

In this section, we will run a relative MMP, and study birational models and divisorial contractions in this MMP.
Let $E$ be the exceptional divisor of a divisorial contraction. To
apply induction on a general fibre of $E$ over its centre, we need
the boundary and its root orders on this fibre. 

\subsection{The birational models}

We use a $\QQ$-factorization obtained by a relative MMP.

\begin{lemma}\label{lem:mmp} (\cite[Corollary~1.37]{Kol13})
Let $(X,\Gamma)$ be a projective dlt pair. There is a projective log
resolution $\pi_0:Y_0\to X$ and a relative MMP over $X$
\[
 (Y_0,\Gamma_0)\dashrightarrow\cdots\dashrightarrow(Y_N,\Gamma_N),
 \qquad \Gamma_0=(\pi_0)_*^{-1}\Gamma+\sum_a E_a ,
\]
where $E_a$ are the exceptional prime divisors of $\pi_0$. Each $(Y_j,\Gamma_j)$ is
$\QQ$-factorial dlt, and every $E_a$ is contracted. The morphism
$\pi_N:Y_N\to X$ is small and crepant and is an isomorphism over
every $\QQ$-factorial open subset of $X$.
\end{lemma}

On a singular model $Y$, an exceptional prime divisor $E$ is in general only
$\QQ$-Cartier. We use the following stack to make the reflexive
pullback of $\OO_Y(E)$ invertible.

\begin{lemma}\label{lem:index}
Let $A_1,\ldots,A_k$ be integral Weil divisors on a normal
$\QQ$-factorial variety $Y$, and let $U\subset Y$ be their common
Cartier locus. There is a normal separated Deligne--Mumford stack
$p_{\rm ind}:\mathfrak Y\to Y$ with coarse space $Y$ and
$\mathfrak Y|_U=U$, such that the reflexivization of each $p_{\rm ind}^*\OO_Y(A_i)$ is invertible. Locally $\mathfrak Y=[Z_0/H]$ with $H$ finite
abelian and $\nu_0:Z_0\to Y$ finite and quasi-\'etale. Assume that
$Y$ is klt. Then $Z_0$ is klt, and each reduced effective $A_i$ pulls
back to a reduced Cartier divisor on $\mathfrak Y$. If $A_i=K_Y$
for some $i$, then $Z_0$ is Gorenstein.
\end{lemma}
\begin{proof}
For each $i$, the sheaf $\OO_Y(A_i)$ is reflexive of rank one and
has an invertible reflexive power, since $Y$ is $\QQ$-factorial. Let $\mathfrak Y_i$ be its index
stack \cite[Definition~5.3.1 and Proposition~5.3.2]{AH09}. Over $U$
each $\mathfrak Y_i$ is isomorphic to $U$. Set
\[
 \mathfrak P=\mathfrak Y_1\times_Y\cdots\times_Y\mathfrak Y_k,
 \qquad \mathfrak Y=(\overline U_{\rm red})^\nu ,
\]
where $\overline U$ is the closure of $U$ in $\mathfrak P$ and
$(\,\cdot\,)^\nu$ denotes normalization.

Take a local trivialization of each $\OO_Y(q_iA_i)$, where $q_i$ is
the Cartier index of $A_i$. The index covers give charts
$\mathfrak Y_i=[Z_i/\mu_{q_i}]$ \cite[Definition~5.19]{KM98}. Over
$U$ the product of the covers $Z_i\to Y$ is a torsor under
$H=\prod_i\mu_{q_i}$. Let $\nu_0:Z_0\to Y$ be the normalization of
$Y$ in this finite \'etale cover of $U$. Then $\mathfrak Y=[Z_0/H]$
on this chart and $Z_0/H=Y$ by normality. The map $\nu_0$ is finite
and is \'etale over $U$. Since $Y\setminus U$ has codimension at least
two, $\nu_0$ is quasi-\'etale.

The tautological line bundle from $\mathfrak Y_i$ agrees over $U$ with
$\OO_U(A_i)$. By normality it is the reflexive pullback of $\OO_Y(A_i)$
on $\mathfrak Y$. If $A_i$ is effective, its section extends across the
complement of $U$. Assume that $Y$ is klt. Since $K_{Z_0}=\nu_0^*K_Y$,
Proposition~\ref{prop:log discrepancies under finite covers} on each
component shows that $Z_0$ is klt. In particular $Z_0$ is Cohen--Macaulay
\cite[Theorem~5.22]{KM98}. For $A_i$ reduced and effective, its
Cartier pullback is generically reduced by quasi-\'etaleness and is
Cohen--Macaulay. Hence it is reduced. If $A_i=K_Y$, the corresponding
line agrees with $\omega_{Z_0}$ in codimension one. By normality
$\omega_{Z_0}$ is invertible, so $Z_0$ is Gorenstein.
\end{proof}

We apply Lemma~\ref{lem:index} to $K_Y$ and the boundary prime divisors.

\begin{lemma}\label{lem:rooted-dlt}
Let $(Y,\Delta+E+D)$ be a $\QQ$-factorial dlt pair with standard
coefficients, where $\lfloor\Delta+E+D\rfloor=E+D$, $E$ is prime,
and $D=\sum_jD_j$. On the index stack for $K_Y$ and the boundary
prime divisors, take a root of order $r$ along each component of $\Delta$ of
coefficient $1-1/r$, and of order $m_j$ along $D_j$. Take no root
along $E$. Let $p:\cS^0\to Y$ be the coarse map and write
$\cE^0,\cD_j^0$ for the reduced inverse images of $E,D_j$.
Set $\cD^0=\sum_j\cD_j^0$. Then the following hold.

\noindent (1) The stack $\cS^0$ has normal Gorenstein klt charts,
and $\cE^0,\cD_j^0$ are Cartier divisors.

\noindent (2) The pair $(\cS^0,\cE^0+\cD^0)$ is dlt on these charts and
\begin{equation}\label{eq:stack-crepant}
 K_{\cS^0}+\cE^0+\cD^0=p^*(K_Y+\Delta+E+D).
\end{equation}

\noindent (3) The pair $(\cS^0,\cE^0)$ is plt near $\cE^0$.
The divisor $\cE^0$ is normal and Gorenstein with trivial generic stabilizer.
\end{lemma}
\begin{proof}
We work locally on $Y$. Lemma~\ref{lem:index} gives a Gorenstein index
cover $\Spec R\to Y$ on which the boundary prime divisors are reduced and Cartier.
If $a_i=0$ are the divisors to be rooted and $s_i$ their orders, the
root chart is
\[
 Z=\Spec R[z_1,\ldots,z_l]/(z_1^{s_1}-a_1,\ldots,z_l^{s_l}-a_l).
\]
It is a complete intersection over $R$, hence Gorenstein. In codimension
one the chart map is either \'etale or, up to an \'etale map, given by
taking a root of a uniformizer. Thus $Z$
is normal by Serre's criterion. The equation of $E$ and the root
coordinates over $D_j$ define generically reduced Cartier divisors.
They are Cohen--Macaulay, hence reduced.

Let $\nu:Z\to Y$ be the chart map, let $E_Z,D_Z$ be the reduced
inverse images of $E,D$, and set $\Lambda=\Delta+\sum_j(1-1/m_j)D_j$.
Let $U$ be the SNC locus of $(Y,\Delta+E+D)$. The index stack is
trivial over $U$ by Lemma~\ref{lem:index}, so
\eqref{eq:root-canonical} with root order one along $E$ gives
\[
 K_Z=\nu^*(K_Y+\Lambda),\qquad
 K_Z+E_Z+D_Z=\nu^*(K_Y+\Delta+E+D).
\]
Both equalities hold on $Z$ since $\nu^{-1}U$ contains every
codimension-one point of $Z$. The equality for $K_Z+E_Z+D_Z$ gives
\eqref{eq:stack-crepant}.

$(Y,\Lambda+E)$ is plt, and $(Y,\Lambda)$ is klt. Since $\nu^*E=E_Z$,
Proposition~\ref{prop:log discrepancies under finite covers} gives that
$Z$ is klt, $(Z,E_Z)$ is plt, and $(Z,E_Z+D_Z)$ is lc. For any lc
centre $W$ of $(Z,E_Z+D_Z)$, the image $\nu(W)$ is an lc centre of
$(Y,\Delta+E+D)$. Thus $\eta_{\nu(W)}\in U$ and $\eta_W\in\nu^{-1}U$.
Over $U$ the chart takes separate
roots of coordinate parameters. Hence $(Z,E_Z+D_Z)$ is SNC at the
generic point of every lc centre and is therefore dlt.

By Proposition~\ref{prop:adjunction} the components of $E_Z$ are
normal. They are disjoint, since an intersection would give a
codimension-two lc centre of the plt pair $(Z,E_Z)$. Thus $E_Z$ is
normal. Moreover, since $E_Z$ is Cartier in $Z$, it is Gorenstein. The index stack and all roots are trivial at $\eta_E$.
\end{proof}

On each model of Lemma~\ref{lem:mmp} take the index stack for $K_{Y_j}$
and the boundary prime divisors, then take their prescribed roots. Write $p_j:\mathcal Y_j\to Y_j$
for the coarse map and set
$D_j=(p_j^{-1}\lfloor\Gamma_j\rfloor)_{\rm red}$.
For a flip let $\theta_j:O_j\xrightarrow{\sim}O'_j$ be the isomorphism
away from the flipping and flipped loci. For a divisorial contraction
$\phi_j$ with exceptional divisor $E$, take
\[
 O_j=Y_j\setminus E,\qquad O'_j=Y_{j+1}\setminus\phi_j(E),\qquad
 \theta_j=\phi_j|_{O_j}.
\]
For a boundary prime divisor $A$ on $O_j$, set $A'=\theta_j(A)$.
Then $\theta_j^*\omega_{O'_j}\simeq\omega_{O_j}$ and
$\theta_j^*(\OO_{O'_j}(A'),1)\simeq(\OO_{O_j}(A),1)$.
A boundary prime divisor on $Y_j$ disjoint from $O_j$ restricts to
the zero divisor on $O_j$, so its index stack and its root are
trivial there. The same holds on $O'_j$ for boundary prime divisors on
$Y_{j+1}$ disjoint from $O'_j$.
Thus $\theta_j$ identifies the normalized products of index stacks
in Lemma~\ref{lem:index}. Taking the same root order on $A$ and $A'$
in Definition~\ref{def:root-stack} gives an isomorphism of pairs over $X$,
\begin{equation}\label{eq:transport}
 (p_j^{-1}O_j,D_j|_{p_j^{-1}O_j})
 \ \xrightarrow{\ \sim\ }\
 (p_{j+1}^{-1}O'_j,D_{j+1}|_{p_{j+1}^{-1}O'_j}).
\end{equation}

\subsection{Adjunction on an exceptional divisor}\label{sec:fibres}

Let $\phi:(Y,\Delta+E+D)\to Y'$ be a divisorial contraction in
Lemma~\ref{lem:mmp}, with exceptional divisor $E$. Take the stack
$p:\cS^0\to Y$ of Lemma~\ref{lem:rooted-dlt} and write
$p_E:\cE^0\to E$ for its restriction. The induced morphism from the
coarse space of $\cE^0$ to $E$ is finite and birational, hence an
isomorphism since $E$ is normal by Proposition~\ref{prop:adjunction}.
We determine the root order on $E$
at each prime divisor $P$.

\begin{lemma}\label{lem:local-ramification}
Set $\Gamma_E=\Diff_E(\Delta+D)$. For a prime divisor $P$ on $E$,
write $\gamma_P=\operatorname{coeff}_P\Gamma_E$ and let $e_P$ be the
ramification index of a finite chart of $p_E$ over $\eta_P$.
The stack $\cS^0$ is smooth and $\cE^0+\cD^0$ is SNC at every
point of $p^{-1}(\eta_P)$.
If $P\not\subset\Supp(D|_E)$, then
\begin{equation}\label{eq:direct-ramification}
 \gamma_P=1-\frac{1}{e_P}.
\end{equation}
If $P\subset\Supp(D|_E)$, there is a unique $j$ with $P\subset D_j|_E$,
and $\gamma_P=1$, $e_P=m_j$. In particular $\Gamma_E$ has standard
coefficients and $\lfloor\Gamma_E\rfloor=D|_E$ is reduced and $\QQ$-Cartier.
\end{lemma}
\begin{proof}
Take a chart $\nu:Z\to Y$ of Lemma~\ref{lem:rooted-dlt} and let $S$
be a component of $E_Z$. Set $h=\nu|_S:S\to E$ and $\Theta=\Diff_S(D_Z)$.
Take a prime divisor $Q$ on $S$ with $h(Q)=P$. Applying the residue isomorphism
of Proposition~\ref{prop:adjunction} to \eqref{eq:stack-crepant}
therefore gives $K_S+\Theta=h^*(K_E+\Gamma_E)$. Proposition~\ref{prop:log discrepancies under finite covers}
applied to $h$ at $Q$ gives
\[
 1-\operatorname{coeff}_Q\Theta=e_P(1-\gamma_P).
\]
If $P\not\subset\Supp(D|_E)$, then $\operatorname{coeff}_Q\Theta=0$,
which proves \eqref{eq:direct-ramification}.
If $P\subset\Supp(D|_E)$, then $P$ is a component of $E\cap D_j$ for
some $j$. By Proposition~\ref{prop:adjunction}, $(Y,\Delta+E+D)$ is
SNC at $\eta_P$. Since $\operatorname{codim}_Y P=2$, the only boundary
components through $\eta_P$ are $E=(x=0)$ and $D_j=(u=0)$ in regular
parameters $x,u$. Adjunction on $x=0$ gives $\gamma_P=1$,
and the root chart $u=v^{m_j}$ is smooth with $e_P=m_j$ on $x=0$.
Thus $\lfloor\Gamma_E\rfloor=D|_E$ is reduced and is $\QQ$-Cartier
because $D$ is $\QQ$-Cartier.
\end{proof}

Let $P_1,\ldots,P_t$ be the components of $\Gamma_E$. By
Lemma~\ref{lem:local-ramification}, after removing a closed subset
of codimension at least two from $E$, the remaining open set $E^\circ$
is smooth, the $P_i^\circ=P_i\cap E^\circ$ are smooth and disjoint,
and $(\cS^0,\cE^0+\cD^0)$ is SNC along $\cE^\circ=p_E^{-1}E^\circ$.
The smooth tame stack $\cE^\circ$ has trivial generic stabilizer.
The root stack with these ramification indices is already smooth,
so \cite[Theorem~1 and Remark~3]{GS15} gives
\begin{equation}\label{eq:boundary-roots}
 \cE^\circ\simeq
 \sqrt[e_{P_1}]{(E^\circ,P_1^\circ)}\times_{E^\circ}\cdots
 \times_{E^\circ}\sqrt[e_{P_t}]{(E^\circ,P_t^\circ)} .
\end{equation}
Here $\cD^0|_{\cE^\circ}$ is the reduced inverse image of $D|_{E^\circ}$.

\subsection{Curves on a general fibre}

Write $E\xrightarrow{s}W_0\to\phi(E)$ for the Stein factorization.
Set $F_w=s^{-1}(w)$ and $F_w^\circ=F_w\cap E^\circ$. For general $w$,
restrict the horizontal components of $\Gamma_E$ to $F_w^\circ$
with their coefficients. Denote the closure of this divisor in $F_w$
by $\Gamma_w$.

\begin{lemma}\label{lem:coarse-fibre}
The divisor $-E$ is $\phi$-ample. There is a nonempty smooth open
$W\subset W_0$ such that for every closed $w\in W$ the following hold.

\noindent (1) The fibre $F_w$ is normal and integral,
$1\leq\dim F_w<\dim Y$, and $F_w\setminus F_w^\circ$ has
codimension at least two in $F_w$.

\noindent (2) The pair $(F_w,\Gamma_w)$ is dlt with standard coefficients, and
\begin{equation}\label{eq:fano-fibre}
 -(K_{F_w}+\Gamma_w)\text{ is ample},\qquad
 B_w=\lfloor\Gamma_w\rfloor=(D|_E)|_{F_w}\text{ is }\QQ\text{-Cartier}.
\end{equation}
\end{lemma}
\begin{proof}
The negativity lemma \cite[Lemma~3.39]{KM98} and $\rho(Y/Y')=1$ give
that $-E$ is $\phi$-ample. Since $E$ is normal and
$s_*\OO_E=\OO_{W_0}$, $K(W_0)$ is algebraically closed in $K(E)$.
The generic fibre of $s$ is geometrically normal and geometrically integral. By
generic flatness and \cite[Theorem~12.2.4(iv), (viii)]{Gro66}, the
general fibre is normal and integral of dimension $\dim E-\dim W_0\geq1$.
After shrinking to a smooth open $W$, the fibre dimension theorem gives
$\operatorname{codim}_{F_w}(F_w\setminus F_w^\circ)\geq2$, and generic
smoothness makes $E^\circ$ and its nonempty boundary strata smooth over $W$. This proves (1).

By Proposition~\ref{prop:adjunction}, $(E,\Gamma_E)$ is dlt. Take a
log resolution $\mu:\widetilde E\to E$ on which every $\mu$-exceptional
divisor has positive log discrepancy, and write
$K_{\widetilde E}+\widetilde\Gamma=\mu^*(K_E+\Gamma_E)$. Shrink $W$
so that $\widetilde E$ and all strata of the boundary and exceptional
divisors are smooth over $W$. Set $\widetilde F_w=(s\mu)^{-1}(w)$.
Then $\mu_w=\mu|_{\widetilde F_w}:\widetilde F_w\to F_w$ is a log
resolution. On $F_w^\circ$, adjunction to the fibre
identifies $K_{F_w}+\Gamma_w$ with $(K_E+\Gamma_E)|_{F_w}$. Since
$F_w$ is normal and $\operatorname{codim}_{F_w}(F_w\setminus F_w^\circ)\geq2$,
this gives
\[
 K_{F_w}+\Gamma_w\sim_\QQ (K_E+\Gamma_E)|_{F_w}
       \sim_\QQ (K_Y+\Delta+E+D)|_{F_w}.
\]
Restricting the discrepancy equality for $\mu$ to $\widetilde F_w$
preserves the coefficients of $\widetilde\Gamma$ on the exceptional
divisors, all strictly less than one. Hence $(F_w,\Gamma_w)$ is dlt.
Since $-(K_Y+\Delta+E+D)$ is $\phi$-ample, $-(K_{F_w}+\Gamma_w)$
is ample. The horizontal components of $\Gamma_E$ restrict on $F_w^\circ$
to disjoint reduced divisors with the same coefficients. Hence
$B_w=(D|_E)|_{F_w}$, which is $\QQ$-Cartier. This proves (2).
\end{proof}

By \eqref{eq:boundary-roots}, the stack $\cE^\circ\times_{E^\circ}F_w^\circ$
has root order $r$ at a component of $\Gamma_w$ of coefficient $1-1/r$
and order $m_j$ at a component of $D_j|_{F_w}$. We use these root
orders in $Q_{\dim F_w}$.

\begin{lemma}\label{lem:fibre}
Assume that $Q_k$ holds for $1\leq k<\dim Y$. There exist an integer
$d>0$ and a nonempty open $\Omega^0\subset\cE^0$ such that the following
property holds. For every $y\in \Omega^0$ there is an
ordinary map $h:(\PP^1,v)\to\cE^0$ with $h(v)\simeq y$. Its coarse
map is nonconstant, its image is contained in a fibre of $sp_E$, and
\begin{equation}\label{eq:d}
 h^*T_{\cE^0}(-\log(\cD^0|_{\cE^0}))\text{ is globally generated},
 \qquad \deg h^*\OO_{\cS^0}(\cE^0)=-d .
\end{equation}
$\cS^0$ is smooth and $\cE^0+\cD^0$ is SNC, along $h(\PP^1)$.
\end{lemma}
\begin{proof}
Take $W$ affine in Lemma~\ref{lem:coarse-fibre} and set
$V=\cE^\circ\cap(sp_E)^{-1}W$, $D_V=\cD^0|_V$. By
\eqref{eq:boundary-roots}, $(V,D_V)$ is smooth with relative SNC
boundary over $W$. For a closed $w\in W$, set $V_w=V\times_W\{w\}$.
The rooted smooth SNC locus
of $(F_w,\Gamma_w)$ contains $V_w$ with complement of codimension
at least two. By \eqref{eq:fano-fibre} and $Q_{\dim F_w}$ there is
a logarithmically free $h:\PP^1\to(V_w,D_V|_{V_w})$.
Set $L=h^*T_{V_w}(-\log(D_V|_{V_w}))$. The exact sequence
\[
 0\longrightarrow L\longrightarrow h^*T_V(-\log D_V)
   \longrightarrow T_wW\otimes\OO_{\PP^1}\longrightarrow0
\]
splits since $H^1(L)=0$ by Lemma~\ref{lem:p1}. Thus $h$ is logarithmically free in $V$.

For $z\notin h^{-1}D_V$, evaluation at $z$ is smooth at $h$ by
Lemma~\ref{lem:free} and Proposition~\ref{prop:sections}\textup{(2)}.
Take a connected neighbourhood $T$ of $h$ in a parameter atlas on
which $\operatorname{ev}_z$ is smooth and every map $h_t$ is
logarithmically free. Set $\Omega^0=\operatorname{ev}_z(T)\subset V$.
This is a nonempty open subset. For $t\in T$, the composite
$sp_E h_t:\PP^1\to W$ is constant since $W$ is affine. Take $k>0$
with $kE$ Cartier. As $p^*(kE)=k\cE^0$,
\[
 -d:=\deg h^*N_{\cE^0/\cS^0}
     =\tfrac1k\deg(p_Eh)^*(\OO_Y(kE)|_E)<0
\]
by $\phi$-ampleness of $-E$ and nonconstancy of $p_Eh$. The degree
$\deg h_t^*N_{\cE^0/\cS^0}$ is constant on $T$, so every $p_Eh_t$
is nonconstant. By the definition of $E^\circ$,
$(\cS^0,\cE^0+\cD^0)$ is SNC along $V\subset\cE^\circ$.
\end{proof}

\section{Smoothing and proof of the main theorem}\label{sec:final}

In this section, we give a comb argument and prove the inductive statement.

\subsection{Smoothing a curve off a divisor}\label{sec:cancellation}

Let $C=C_1\cup C_2$ be the union of two curves at smooth marked
points identified to a node $p$. Its local ring at $p$ is
$\OO_{C_1,p}\times_{\kappa(p)}\OO_{C_2,p}$. By the pushout property
\cite[Lemma~4.5]{AHHLR24}, a map from $C$ to a
separated Deligne--Mumford stack is given by maps on $C_1,C_2$ and
an isomorphism between their values at $p$.

\begin{proposition}\label{prop:cancellation}
Let $\cS$ be a smooth separated Deligne--Mumford stack with
quasi-projective coarse space. Let $E+D_1+\cdots + D_s$ be a reduced SNC divisor on $\cS$. Set $D'=\sum_jD_j$. Assume that the following hold.

\noindent (1) $f:C_0\simeq\PP^1\to(\cS,E+D')$ is logarithmically free.
We have $f^*E=\sum_{a=1}^r[p_a]$ such that $p_a$ are distinct and
$f(p_a)\notin D'$.

\noindent (2) For every $a$ there is an ordinary map $h_a:(C_a,v_a)\to E$
from a marked $\PP^1$ and an isomorphism $h_a(v_a)\simeq f(p_a)$.

\noindent (3) For every $a$, $h_a(C_a)\not\subset D'|_E$,
$h_a^*T_E(-\log(D'|_E))$ is globally generated, and
$h_a^*N_{E/\cS}\simeq\OO_{C_a}(-1)$.

Attach $C_a$ to $C_0$ by identifying $v_a$ with $p_a$.
The resulting map deforms to a logarithmically free ordinary map
$\PP^1\to(\cS\setminus E,D'|_{\cS\setminus E})$.
\end{proposition}
\begin{proof}
For $r=0$ take $f$. Otherwise let $u:C\to\cS$ be the glued map and
take $z\in C_0$ with $f(z)\notin E\cup D'$. Set $L=T_{\cS}(-\log D')$,
$L_a=L|_{f(p_a)}$, $N_a=N_{E/\cS}|_{f(p_a)}$, and
$i_a:\{p_a\}\hookrightarrow C_0$. Then we have
\begin{align}
 0&\longrightarrow f^*T_{\cS}(-\log(E+D'))\longrightarrow f^*L
   \longrightarrow\bigoplus_a(i_a)_*N_a\longrightarrow0,
       \label{eq:contact-normal}\\
 0&\longrightarrow h_a^*T_E(-\log(D'|_E))\longrightarrow h_a^*L
   \longrightarrow\OO_{C_a}(-1)\longrightarrow0.
       \label{eq:tail-normal}
\end{align}
Here \eqref{eq:tail-normal} holds by $h_a^*N_{E/\cS}\simeq \OO_{C_a}(-1)$.

For each $a$, let $q_a:L_a\to N_a$ be the normal projection. By
\textup{(1)} and Lemma~\ref{lem:p1} we have
\[
 H^1\bigl(C_0,f^*T_{\cS}(-\log(E+D'))(-z)\bigr)=0 .
\]
Since $z\ne p_a$ for all $a$, twisting \eqref{eq:contact-normal} by $-z$ gives
$H^1(C_0,f^*L(-z))=0$ and a surjection
\[
 H^0(C_0,f^*L(-z))\longrightarrow\bigoplus_a N_a,\qquad
 s\longmapsto \bigl(q_a(s(p_a))\bigr)_a .
\]
On $C_a$, set $F_a=h_a^*T_E(-\log(D'|_E))$. By \textup{(3)}, $F_a$
is globally generated, so $H^1(C_a,F_a)=0$ by Lemma~\ref{lem:p1}.
Since $H^1(C_a,\OO_{C_a}(-1))=0$, \eqref{eq:tail-normal} gives
$H^1(C_a,h_a^*L)=0$.

The normalization sequence for $C$ now identifies $H^1(C,u^*L(-z))$
with the cokernel of
\[
 \begin{aligned}
 \delta:\ H^0(C_0,f^*L(-z))\oplus\bigoplus_a H^0(C_a,h_a^*L)
       &\longrightarrow \bigoplus_a L_a,\\
 (s,(t_a)_a)&\longmapsto \bigl(s(p_a)-t_a(v_a)\bigr)_a .
 \end{aligned}
\]
Given $(w_a)_a\in\bigoplus_a L_a$, take $s$ with
$q_a(s(p_a))=q_a(w_a)$ for every $a$. Then
$s(p_a)-w_a\in\ker q_a$. Since $f(p_a)\notin D'$, \textup{(2)} gives $\ker q_a=F_a|_{v_a}$. Global generation of $F_a$
gives a section $t_a$ with $t_a(v_a)=s(p_a)-w_a$. Via
\eqref{eq:tail-normal} we may view $t_a$ as a section of $h_a^*L$.
Thus $\delta(s,(t_a)_a)=(w_a)_a$, proving $H^1(C,u^*L(-z))=0$.
The map $u^*L(-z)\to u^*T_{\cS}$ is an isomorphism away from
$\{z\}\cup u^{-1}D'$, a finite set of smooth points by \textup{(1)}
and \textup{(3)}. It is therefore injective with finite-length
cokernel, so $H^1(C,u^*T_{\cS})=0$.

Let $\cC\to\Aff^1$ be the family obtained by blowing up
$C_0\times\Aff^1$ at the points $(p_a,0)$. Identify the exceptional
line at $p_a$ with $C_a$, taking its intersection with the strict
transform of $C_0\times\{0\}$ to $v_a$. The central fibre is then $C$,
and the strict transform of $\{z\}\times\Aff^1$ gives a section $\sigma$.
By $H^1(C,u^*T_{\cS})=0$ and Proposition~\ref{prop:sections}\textup{(1)},
$u$ extends to $G:\cC_B\to\cS$ over an \'etale neighbourhood
$(B,b_0)\to(\Aff^1,0)$. Denote the induced section on $B$ by $\sigma_B$.
For a general $t\in B$, let $g$ be the restriction of $G$ to the fibre
$C_t$ and set $z_t=\sigma_B(t)$. Since $u(z)\notin E\cup D'$, the
open set $(G\sigma_B)^{-1}(\cS\setminus(E\cup D'))$ contains $b_0$.
Thus $g(z_t)\notin E\cup D'$. Upper semicontinuity
for $G^*L(-\sigma_B(B))$ and $H^1(C,u^*L(-z))=0$ give
$H^1(C_t,g^*L(-z_t))=0$.
As $C_t\simeq\PP^1$, Lemma~\ref{lem:p1} gives global generation of $g^*L$.

Note that
\begin{equation}\label{eq:zero-intersection}
 \deg g^*E=\deg f^*E+\sum_a\deg h_a^*N_{E/\cS}=r-r=0.
\end{equation}
The effective divisor $g^*E$ is therefore zero. Let $\pi:\cS\to S$ be
the coarse map and $H$ an ample line bundle on $S$. Then
$\deg(\pi g)^*H=\deg(\pi f)^*H+\sum_a\deg(\pi h_a)^*H
\geq\deg(\pi f)^*H>0$. Thus $\pi g$ is nonconstant and $g$ is
logarithmically free in $(\cS\setminus E,D'|_{\cS\setminus E})$.
\end{proof}

\subsection{Root orders and the induction on dimension}

Fix the MMP $Y_0\dashrightarrow\cdots\dashrightarrow Y_N$ of
Lemma~\ref{lem:mmp} for an $n$-dimensional pair $(X,\Gamma=\Delta+B)$ as in
Theorem~\ref{thm:induction}. Write $\phi_{j_i}:Y_{j_i}\to Y_{j_i+1}$,
$j_1<\cdots<j_q$, for the divisorial contractions, and let
$E_i\subset Y_{j_i}$ be the exceptional prime divisor of $\phi_{j_i}$.
Use the same notation for its strict transform on each earlier model.
On each $Y_j$ take the index stack for $K_{Y_j}$ and the prime divisors
in $\Supp\Gamma_j$ as in Lemma~\ref{lem:index}. Along the Cartier
pullbacks of the strict transforms of $\Delta_k$, $B_\ell$ and $E_i$,
take roots of orders $r_k$, $m_\ell$ and $M_i$, respectively.
Write $p_j:\mathcal Y_j\to Y_j$ for the coarse map of the stack and $D_j=(p_j^{-1}\lfloor\Gamma_j\rfloor)_{\rm red}$ as in
\eqref{eq:transport}. The orders $r_k,m_\ell$ come from
\eqref{eq:root-open}. We determine the orders $M_i$ as follows.

\begin{proposition}\label{prop:root-orders}
Assume that $Q_k$ holds for $1\leq k<n$. For every choice of
$m_1,\ldots,m_b$, there exist positive integers $M_1,\ldots,M_q$
such that the following property holds. For each $i$, set
\[
 \cE_i=(p_{j_i}^{-1}E_i)_{\rm red},\qquad \cD_i=D_{j_i}-\cE_i .
\]
There is a nonempty open $\Omega_i\subset\cE_i$ such that for every
geometric point $y$ of $\Omega_i$ there is an ordinary map
$h:(\PP^1,v)\to\cE_i$ with $h(v)\simeq y$ and $\phi_{j_i}p_{j_i}h$
constant. Along $h(\PP^1)$ the stack $\mathcal Y_{j_i}$ is smooth
and $\cE_i+\cD_i$ is SNC. Moreover, $h^*N_{\cE_i/\mathcal Y_{j_i}}\simeq\OO_{\PP^1}(-1)$ and $h^*T_{\cE_i}(-\log(\cD_i|_{\cE_i}))$ is globally generated.
\end{proposition}
\begin{proof}
We determine $M_q,\ldots,M_1$ in this order. On $Y_{j_i}$ the
exceptional prime divisors over $X$ are $E_i,\ldots,E_q$. To determine the order for $E_i$, we can assume that the orders on $E_{i+1},\ldots,E_q$ are already fixed.
Let $p_i^0:\cS_i^0\to Y_{j_i}$ be the stack of Lemma~\ref{lem:rooted-dlt}
with orders $r_k,m_\ell$ on the transforms of $\Delta_k,B_\ell$,
the fixed orders on $E_{i+1},\ldots,E_q$, and order one on $E_i$.
Write $\cE_i^0=((p_i^0)^{-1}E_i)_{\rm red}$ and let $\cD_i^0$ be the
reduced inverse image of $\lfloor\Gamma_{j_i}\rfloor-E_i$.

Lemma~\ref{lem:fibre} gives a nonempty open $\Omega_i^0\subset\cE_i^0$
and an integer $d_i>0$. Through each point of $\Omega_i^0$ the lemma
gives a marked curve on $\cE_i^0$ of normal degree $-d_i$.
Set
\begin{equation}\label{eq:root-choice}
 M_i=d_i,\qquad
 \lambda_i:\mathcal Y_{j_i}=\sqrt[d_i]{(\cS_i^0,\cE_i^0)}
                 \longrightarrow\cS_i^0 .
\end{equation}
The zero divisor of the universal root section is $\cE_i$ and
$p_{j_i}=p_i^0\lambda_i$. Set
$\Omega_i=(\lambda_i|_{\cE_i})^{-1}(\Omega_i^0)$.

Fix a geometric point $y$ of $\Omega_i$ and take such a curve
$h^0:(\PP^1,v)\to\cE_i^0$ with $h^0(v)\simeq\lambda_i(y)$.
Since $(h^0)^*N_{\cE_i^0/\cS_i^0}$ has degree $-d_i$, there is
an isomorphism
\begin{equation}\label{eq:normal-root}
 \alpha:\OO_{\PP^1}(-1)^{\otimes d_i}
       \xrightarrow{\ \sim\ }(h^0)^*N_{\cE_i^0/\cS_i^0}.
\end{equation}
By Definition~\ref{def:root-stack},
the triple $(\OO_{\PP^1}(-1),0,\alpha)$ on $\PP^1$ gives a lift
$h:\PP^1\to\cE_i$ of $h^0$. The universal root line is
$\OO_{\mathcal Y_{j_i}}(\cE_i)$, so its restriction to $\cE_i$
is $N_{\cE_i/\mathcal Y_{j_i}}$ and its pullback by $h$ is
$\OO_{\PP^1}(-1)$. Over the geometric point $\lambda_i(y)$, any two
$d_i$-th roots of the same line are isomorphic. Hence $h(v)\simeq y$.
Moreover, $\phi_{j_i}p_{j_i}h=\phi_{j_i}p_i^0h^0$ is constant
by Lemma~\ref{lem:fibre}.

Along $h^0(\PP^1)$ the pair $(\cS_i^0,\cE_i^0+\cD_i^0)$ is SNC
by Lemma~\ref{lem:fibre}. Take \'etale coordinates $x_1,\ldots,x_n$
in which $\cE_i^0=(x_1=0)$ and the components of $\cD_i^0$ through
the point are $x_k=0$ for $2\leq k\leq s$.
The root chart has coordinates $u,x_2,\ldots,x_n$ with $x_1=u^{d_i}$.
It is smooth with SNC boundary $(ux_2\cdots x_s=0)$.
On $u=0$ the group $\mu_{d_i}$ fixes $x_2,\ldots,x_n$, so
$d(\lambda_i|_{\cE_i})$ induces an isomorphism
\[
 h^*T_{\cE_i}(-\log(\cD_i|_{\cE_i}))
 \simeq (h^0)^*T_{\cE_i^0}(-\log(\cD_i^0|_{\cE_i^0})).
\]
The bundle on the right is globally generated by Lemma~\ref{lem:fibre}.
\end{proof}

\begin{proof}[Proof of Theorem~\ref{thm:induction}]
Proposition~\ref{prop:curve} gives $Q_1$. Assume that $n\geq2$ and
that $Q_k$ holds for $1\leq k<n$. Fix the orders $m_\ell$ on $B_\ell$
and take the orders $M_i$ of Proposition~\ref{prop:root-orders}.
For each $\pi_j:Y_j\to X$, set $\rho_j=\pi_jp_j$.
Let $\mathcal V_j\subset\mathcal Y_j$ be the open locus where
$\mathcal Y_j$ is smooth and $D_j$ is SNC. Fix a very ample divisor
$H$ on $X$ and set $H_j=\rho_j^*\OO_X(H)$.

\smallskip\noindent\textbf{Step~1.} In this step, we construct a
logarithmically free curve $f_0$ on $(\mathcal Y_0,D_0)$ with
$\deg f_0^*H_0>0$. Let $X^\circ$ be the intersection of $U$ with
the isomorphism locus of $\pi_0$ in $X$. Its complement has
codimension at least two. For a general $x\in X^\circ\setminus\Supp\Gamma$, take
the intersection $C$ of $n-1$ general members of $|H|$ through $x$.
Since $\dim(X\setminus X^\circ)\leq n-2$, $C\subset X^\circ$,
and $C$ is smooth and integral by Bertini's theorem. Its strict
transform $\widetilde C\subset Y_0$ maps isomorphically to $C$.
Since $Y_0$ is smooth and $\Supp\Gamma_0$ is SNC,
\eqref{eq:root-canonical} gives
\begin{equation}\label{eq:initial-canonical}
 K_{\mathcal Y_0}+D_0=p_0^*(K_{Y_0}+\Gamma_0).
\end{equation}
By Lemma~\ref{lem:ct-lift} there are a finite cover
$\tau:C'\to\widetilde C$ from a smooth connected projective curve and
a map $b:C'\to\mathcal Y_0$ lifting
the composite $C'\to\widetilde C\hookrightarrow Y_0$. Since $\pi_0$
is an isomorphism along $\widetilde C$, \eqref{eq:initial-canonical}
gives
\[
 \deg b^*(K_{\mathcal Y_0}+D_0)
   =(\deg\tau)(K_X+\Gamma)\cdot H^{n-1}<0 .
\]
The inequality follows from bigness of $-(K_X+\Gamma)$. Since
$p_0b(C')=\widetilde C$ contains the unique point over $x$,
Lemma~\ref{lem:initial} gives a logarithmically free map
$f_0:\PP^1\to(\mathcal Y_0,D_0)$. By Lemma~\ref{lem:free} and
Proposition~\ref{prop:sections}\textup{(2)}, evaluation at $0$ is
smooth at $f_0$ on the open locus of logarithmically free maps.
After deforming $f_0$, we may therefore take $p_0f_0(0)\in\pi_0^{-1}X^\circ$.
If $\rho_0f_0$ were constant, the image of $p_0f_0$ would lie in the
single-point fibre $\pi_0^{-1}(\rho_0f_0(0))$. Since $p_0f_0$ is
nonconstant, we obtain $\deg f_0^*H_0>0$.

\smallskip\noindent\textbf{Step~2.} In this step, we pass a
logarithmically free curve of positive degree over $X$ from $\mathcal V_j$
to $\mathcal V_{j+1}$. Let $f_j:\PP^1\to(\mathcal V_j,D_j|_{\mathcal V_j})$
be such a curve. If $Y_j\dashrightarrow Y_{j+1}$ is a flip, the
flipping locus $Y_j\setminus O_j$ has codimension at least two.
Its inverse image in $\mathcal V_j$ has the same codimension bound
since $p_j$ has finite charts. By Lemma~\ref{lem:free}\textup{(2)},
a logarithmically free deformation $g_j$ of $f_j$ lies in $p_j^{-1}O_j$.
In this family $\deg g_j^*H_j=\deg f_j^*H_j$.

For a divisorial contraction, write $j=j_i$ and $D_j=\cE_i+\cD_i$.
Since $\cE_i$ is irreducible, both $\cE_i\setminus\Omega_i$ and
$\cE_i\cap\cD_i$ are proper closed subsets of the divisor $\cE_i$.
They have codimension at least two in $\mathcal Y_j$. By
Lemma~\ref{lem:free}\textup{(2)}, we may take a logarithmically free
deformation $\widetilde f_j$ of $f_j$ such that
\[
 \widetilde f_j^*\cE_i=\sum_{a=1}^r[p_a],\qquad
 p_a\ne p_b\ (a\ne b),\qquad
 \widetilde f_j(p_a)\in\Omega_i\setminus\cD_i .
\]
At each $p_a$, take the marked curve $h_a:(\PP^1,v_a)\to\cE_i$
of Proposition~\ref{prop:root-orders} such that $h_a(v_a)\simeq\widetilde f_j(p_a)$.
Its image lies in $\mathcal V_j$ and meets $\cE_i\setminus\cD_i$.
Moreover, $h_a^*T_{\cE_i}(-\log(\cD_i|_{\cE_i}))$ is globally generated
and $h_a^*N_{\cE_i/\mathcal Y_j}\simeq\OO_{\PP^1}(-1)$.
Proposition~\ref{prop:cancellation}
applied in $\mathcal V_j$ therefore gives a
logarithmically free $g_j$ in $\mathcal V_j\setminus\cE_i$.
Since $\phi_jp_jh_a$ is constant and $\rho_j=\pi_{j+1}\phi_jp_j$,
each $h_a^*H_j$ has degree zero. Then we have
\[
 \deg g_j^*H_j
   =\deg\widetilde f_j^*H_j+\sum_{a=1}^r\deg h_a^*H_j
   =\deg f_j^*H_j .
\]
This also holds for $r=0$, with $g_j=\widetilde f_j$.
In both cases $g_j$ lies in $p_j^{-1}O_j$. The isomorphism of pairs
\eqref{eq:transport} over $X$ carries $g_j$ to a logarithmically free
map $f_{j+1}$ in $\mathcal V_{j+1}$ and identifies the restrictions
of $H_j$ and $H_{j+1}$. Thus $\deg f_{j+1}^*H_{j+1}=\deg f_j^*H_j>0$.

\smallskip\noindent\textbf{Step~3.} In this step, we obtain the
required curve on $(\cU,B_{\cU})$. By Lemma~\ref{lem:mmp},
$\pi_N$ is small and is an isomorphism over the smooth open set $U$.
All the boundary divisors and $K_U$ are Cartier on $U$, so the index
stack is trivial there by Lemma~\ref{lem:index}. Thus $\rho_N^{-1}U$
is the stack $\cU$ of \eqref{eq:root-open},
with $D_N|_{\cU}=B_{\cU}$. Since $X\setminus U$ has codimension at
least two and $\pi_N$ contracts no divisor, $Y_N\setminus\pi_N^{-1}U$
has codimension at least two. The same holds for
$\mathcal V_N\setminus\cU$ since $p_N$ has finite charts.
By Lemma~\ref{lem:free}\textup{(2)}, a general deformation $f$ of
$f_N$ lies in $\cU$ and is logarithmically free there. Evaluation at
$0$ is smooth at $f$ by Lemma~\ref{lem:free} and
Proposition~\ref{prop:sections}\textup{(2)}, so the ordinary
deformations of $f$ cover a dense open subset of $\cU$. For each
closed $Z\subset\cU$ of codimension at least two,
Lemma~\ref{lem:free}\textup{(2)} applied in $\cU$ gives a general
deformation of $f$ avoiding $Z$.
\end{proof}

\begin{proof}[Proof of Theorem~\ref{mainthm}]
Apply Theorem~\ref{thm:induction} to the klt pair $(X,0)$ of
dimension $n\geq1$. Here $-K_X$ is ample, $U=X_{\reg}$ and
$(\cU,B_{\cU})=(X_{\reg},0)$ by \eqref{eq:root-open}. Thus the
resulting map $f:\PP^1\to X_{\reg}$ is nonconstant and
$f^*T_{X_{\reg}}$ is globally generated.
\end{proof}

\bibliographystyle{amsalpha}
\bibliography{v13}

@misc{AH09,
  author = {Abramovich, Dan and Hassett, Brendan},
  title = {Stable varieties with a twist},
  year = {2009}, howpublished = {arXiv:0904.2797v1},
  url = {https://arxiv.org/abs/0904.2797v1},
  note = {Version 1, 17 April 2009}
}

@article{AHHLR24,
  author = {Alper, Jarod and Hall, Jack and Halpern-Leistner, Daniel and Rydh, David},
  title = {Artin algebraization for pairs with applications to the local structure
           of stacks and {Ferrand} pushouts},
  journal = {Forum Math. Sigma},
  volume = {12}, year = {2024}, pages = {e20, 1--25},
  doi = {10.1017/fms.2023.60},
  url = {https://doi.org/10.1017/fms.2023.60}
}

@article{AV02,
  author = {Abramovich, Dan and Vistoli, Angelo},
  title = {Compactifying the space of stable maps},
  journal = {J. Amer. Math. Soc.},
  year = {2002}, volume = {15}, number = {1}, pages = {27--75},
  doi = {10.1090/S0894-0347-01-00380-0},
  url = {https://arxiv.org/abs/math/9908167}
}

@misc{BCHM08,
  author = {Birkar, Caucher and Cascini, Paolo and Hacon, Christopher D. and McKernan, James},
  title = {Existence of minimal models for varieties of log general type},
  year = {2008}, howpublished = {arXiv:math/0610203v2},
  url = {https://arxiv.org/abs/math/0610203v2},
  note = {Version 2, 14 August 2008}
}

@misc{CT26a,
  author = {Chen, Jiun-Cheng and Tseng, Hsian-Hua},
  title = {On rational curves relative to a divisor},
  year = {2026}, howpublished = {arXiv:2611.05809v1},
  url = {https://arxiv.org/abs/2611.05809v1},
  note = {Version 1, 4 November 2026}
}

@misc{CT26b,
  author = {Chen, Jiun-Cheng and Tseng, Hsian-Hua},
  title = {Smooth weak log {Fano} orbifolds are {Campana} rationally connected},
  year = {2026}, howpublished = {arXiv:2611.15191v1},
  url = {https://arxiv.org/abs/2611.15191v1},
  note = {Version 1, submitted 11 November 2026. Manuscript dated 13 November 2026}
}

@article{Cad07,
  author = {Cadman, Charles},
  title = {Using stacks to impose tangency conditions on curves},
  journal = {Amer. J. Math.},
  year = {2007}, volume = {129}, number = {2}, pages = {405--427},
  doi = {10.1353/ajm.2007.0007},
  url = {https://arxiv.org/abs/math/0312349v3},
  note = {Author version, arXiv:math/0312349v3}
}

@misc{GS15,
  author = {Geraschenko, Anton and Satriano, Matthew},
  title = {A ``bottom up'' characterization of smooth {Deligne--Mumford} stacks},
  year = {2015}, howpublished = {arXiv:1503.05478v1},
  url = {https://arxiv.org/abs/1503.05478v1},
  note = {Version 1, 18 March 2015}
}

@article{Gro66,
  author = {Grothendieck, Alexander},
  title = {{\'E}l{\'e}ments de g{\'e}om{\'e}trie alg{\'e}brique : {IV}.
           {\'E}tude locale des sch{\'e}mas et des morphismes de sch{\'e}mas,
           Troisi{\`e}me partie},
  journal = {Publ. Math. Inst. Hautes {\'E}tudes Sci.},
  volume = {28}, year = {1966}, pages = {5--255},
  doi = {10.1007/BF02684343},
  url = {https://www.numdam.org/item/PMIHES_1966__28__5_0/},
  note = {Written with the collaboration of Jean Dieudonn{\'e}}
}

@misc{GX26,
  author = {Guo, Bin and Xie, Song-Yan},
  title = {Constructing rational curves via jets on projective varieties
           with non-nef canonical bundle},
  year = {2026},
  howpublished = {\href{https://arxiv.org/abs/2609.28465v1}{arXiv:2609.28465v1}},
  note = {Appendix B by Bin Dong, Guoxiong Gao, Zeming Sun, and Bin Wu}
}

@misc{JLR25,
  author = {Jovinelly, Eric and Lehmann, Brian and Riedl, Eric},
  title = {Free curves and fundamental groups},
  year = {2025},
  howpublished = {Preprint},
  note = {\href{https://arxiv.org/abs/2510.27031v1}{arXiv:2510.27031v1}, to appear in Duke Math. J.}
}

@misc{JLR26,
  author = {Jovinelly, Eric and Lehmann, Brian and Riedl, Eric},
  title = {Finding large families of rational curves through {Bend-and-Break}},
  year = {2026},
  howpublished = {\href{https://arxiv.org/abs/2608.27696v1}{arXiv:2608.27696v1}}
}

@misc{KM97,
  author = {Keel, Se{\'a}n and McKernan, James},
  title = {Rational curves on quasi-projective surfaces},
  year = {1997}, howpublished = {Preprint of 25 July},
  url = {https://arxiv.org/abs/alg-geom/9707016v1},
  note = {\href{https://arxiv.org/abs/alg-geom/9707016v1}{arXiv:alg-geom/9707016v1}}
}

@book{KM98,
  author = {Koll{\'a}r, J{\'a}nos and Mori, Shigefumi},
  title = {Birational Geometry of Algebraic Varieties},
  series = {Cambridge Tracts in Mathematics}, volume = {134},
  publisher = {Cambridge University Press}, address = {Cambridge}, year = {1998},
  note = {With the collaboration of C.~H. Clemens and A. Corti},
  doi = {10.1017/CBO9780511662560}
}

@book{Kol13,
  author = {Koll{\'a}r, J{\'a}nos},
  title = {Singularities of the Minimal Model Program},
  series = {Cambridge Tracts in Mathematics}, volume = {200},
  publisher = {Cambridge University Press}, address = {Cambridge}, year = {2013},
  note = {With the collaboration of S{\'a}ndor Kov{\'a}cs},
  doi = {10.1017/CBO9781139547895}
}

@article{Ols06,
  author = {Olsson, Martin C.},
  title = {{Hom}-stacks and restriction of scalars},
  journal = {Duke Math. J.}, volume = {134}, number = {1},
  year = {2006}, pages = {139--164},
  doi = {10.1215/S0012-7094-06-13414-2},
  url = {https://drive.google.com/file/d/1rnWxC5VEh7CrufXjuPcxVjrQz2XeZUSr/view},
  note = {Author copy}
}

@misc{Sta26,
  author = {The {Stacks Project Authors}},
  title = {The {Stacks Project}},
  year = {2026},
  note = {Tags \href{https://stacks.math.columbia.edu/tag/054K}{054K} and
          \href{https://stacks.math.columbia.edu/tag/05XD}{05XD}}
}
\end{document}